\documentclass[psamsfonts]{amsart}

\newtheorem{theorem}{Theorem}[section]
\newtheorem{lemma}[theorem]{Lemma}
\newtheorem{proposition}[theorem]{Proposition}
\newtheorem{corollary}[theorem]{Corollary}

\theoremstyle{definition}
\newtheorem{definition}[theorem]{Definition}
\newtheorem{example}[theorem]{Example}

\newtheorem{remark}[theorem]{Remark}

\theoremstyle{remark}

\numberwithin{equation}{section}

\newcommand{\NN}{\mathbb{N}}
\newcommand{\ZZ}{\mathbb{Z}}
\newcommand{\QQ}{\mathbb{Q}}

\renewcommand{\AA}{\mathbb{A}}

\newcommand  {\foa}     {\mathfrak{a}}

\newcommand  {\fob}     {\mathfrak{b}}

\newcommand  {\fom}     {\mathfrak{m}}

\newcommand  {\foq}     {\mathfrak{q}}
\newcommand  {\fon}     {\mathfrak{n}}

\newcommand  {\Ann}     {\operatorname{Ann}}

\newcommand  {\Hom}     {\operatorname{Hom}}

\renewcommand  {\ker }  {\operatorname{kern}}

\newcommand  {\lra}     {\longrightarrow}

\renewcommand{\O}       {\mathcal{O}}

\newcommand  {\ra}      {\rightarrow}

\newcommand  {\Spec}    {\operatorname{Spec}}

\newcommand {\soclo} {{sc}}
\newcommand {\pasoclo} {\star}
\newcommand {\defop} {\mu}

\def\mydate{\number\day\space\ifcase\month \or January\or February\or March\or April\or May\or
June\or July\or August\or September\or October\or November\or
December\fi \space\number\year}

\usepackage{amscd}
\usepackage{amssymb}

\begin{document}

\title[How to rescue solid closure]
{How to rescue solid closure}

% Remove or comment out any unused author tags.
% author one information

\author[Holger Brenner]{Holger Brenner}
\address{Mathematische Fakult\"at, Ruhr-Universit\"at, 
               44780 Bochum, Germany}
%\curraddr{}
\email{brenner@cobra.ruhr-uni-bochum.de}

%\thanks{}

%\subjclass{14H60, 14J26, 13A35}
%\date{}

% at present the "communicated by" line appears only in ERA and PROC
%\commby{}

\dedicatory{\mydate}

\begin{abstract}
We define a closure operation for ideals in a commutative ring
which has all the good properties of solid closure (at least in the
case of equal characteristic) but such that also
every ideal in a regular ring is closed. This gives in particular
a kind of tight closure theory in characteristic zero
without referring to positive characteristic.
\end{abstract}

\maketitle

%===========================================================
\section*{Introduction}

The theory of tight closure,
introduced by Hochster and Huneke
(\cite{hohutight},\cite{hohubriancon},\cite{hunekeapplication},
\cite{hunekeparameter},\cite{brunsherzogrev}),
is defined for rings containing
a field of positive characteristic making use of the Frobenius
endomorphism.
The main applications of tight closure
are to homological conjectures. For example,
with this technique it is quite easy to proof that the invariant ring
of an action by a linearly
reductive group on a polynomial ring is Cohen-Macaulay.
This rests upon three facts of tight closure: colon-capturing,
persistence and the fact that every ideal in a regular ring is tightly closed.

There are several ways to extend this notion to
rings containing a field of characteristic zero,
by declaring that an element belongs to the tight closure if and only if
this holds for almost all points with residue class field of positive
characteristic, where
the relevant datas have to be expressed
in an algebra over a finitely generated $\ZZ$-algebra.
Most results from tight closure pass then from the case of
positive characteristic to characteristic zero.

The theory of solid closure was introduced by Hochster in \cite{hochstersolid}
in order to give a closure operation defined for every Noetherian ring
without referring to positive characteristic.
It coincides with tight closure if the ring contains a field
of positive characteristic under some mild finiteness conditions,
and it contains tight closure also in the case of
a field of characteristic zero.
However, a computation of Roberts in \cite{robertscomp}
showed that regular rings of dimension $\geq 3$ in characteristic zero
are not solidly closed. This ``discouraging'' result shows that
solid closure
``does not have the right properties in equal characteristic zero''
(\cite[Introduction]{hochstertightsolid}).

This paper proposes another closure operation,
called parasolid closure.
Roughly speaking, an element $h$
belongs to the parasolid closure of $(f_1, \ldots , f_n) \subseteq R$
(where $R$ is local of dimension $d$),
if the forcing algebra $A$ for these datas is (universally) parasolid:
this means that every canonical element (given by parameters in $R$)
coming from $H^d_\fom (R)$ survive
in $H^d_\fom (A)$.

We list the main properties of this closure
operation (compare \cite[Introduction]{hunekeapplication}).
The parasolid closure is persistent and
it lies inside solid closure. Every ideal in a regular ring
is parasolidly closed.

If the ring contains a field,
then the contraction from a finite extension belongs to the closure.
This is not clear in mixed characteristic, and I have to admit
that our closure operation has no effect on the homological conjectures
in mixed characteristic.

Over a field of positive characteristic it coincides with solid closure
and hence also with tight closure under some mild conditions.
Though the definition
of parasolid closure does not refer to positive characteristic,
reduction mod $p$ is anyway a useful method to prove results
for parasolid closure in equal characteristic zero.
The theorem of Brian\c{c}on-Skoda holds for parasolid closure
in equal characteristic zero and we prove this
with the help of the finiteness theorem of
Hochster (\cite[Theorem 8.4.1]{brunsherzogrev}.
Moreover, parasolid closure contains tight closure in characteristic zero,
hence also colon capturing holds for it.

For a complete local Gorenstein ring we also define an ideal
which coincides with the tight closure test ideal in
positive characteristic.

The content of this paper is as follows.
Section \ref{parasolid} gives the concept of a parasolid algebra
and some characterizations in positive characteristic and in low dimensions. 
Section \ref{ordersec} gives criteria for an algebra $A$ over a
complete local Gorenstein ring $R$ to be parasolid
in terms of the module-sections
$\Hom_R(A,R)$ and the order ideal in $R$.
We show that a regular complete local
ring is in every parasolid algebra a direct summand.

Section \ref{parasolidclosure} introduce the closure operation and
establishes its basic properties.
In section \ref{regularsec} we show that every ideal in a regular ring
is closed. Rings with this property will be called pararegular.
Section \ref{solidparasolidsec} is concerned with the relation between
solid and parasolid closure and with tight closure in positive
characteristic.

In section \ref{testandorder} we look at the intersection of the order
ideals of all parasolid forcing algebras. We show that
for a complete local Gorenstein ring over a field of positive characteristic
this intersection is the same as the tight closure test ideal.
In section \ref{pararegularsec} we give characterizations
for pararegular rings in the case of a complete Gorenstein ring.

In section \ref{expresssec} we express the datas which say that a paraclass
vanishes in terms of a finitely generated $\ZZ$-algebra and equations.
This gives relations
between the parasolid closure in positive characteristic,
in zero characteristic and in mixed characteristic.
We show that the equations which define over a field a
quotient singularity in dimension two
yield also in mixed characteristic examples of pararegular,
but non regular rings (this is a kind of reduction to the case of a field).
We obtain the Theorem of Brian\c{c}on-Skoda for parasolid closure in equal
characteristic zero by reduction to positive characteristic
using Hochster's finiteness Theorem. Furthermore
we get from a result of Koh
\cite{kohsuper} a finite type criterion for universally parasolid algebras.

In section \ref{generic} we describe results how
parasolid algebras and the parasolid closure
behave in the fibers of a family, and derive from this as a special case
that the tight closure of an ideal in characteristic zero is contained
inside the parasolid closure.

\section{parasolid algebras}
\label{parasolid}

Let $R$ denote a local Noetherian ring of dimension $d$ with maximal
ideal $\fom$.
Due to a theorem of Grothendieck we know that
the local cohomology $H^d_\fom (R) \neq 0$.
Elements $x_1, \ldots, x_d \in R$ are called parameters
if $V(x_1, \ldots , x_d) =V(\fom)$.
A system of parameters defines via $\rm \check{C}$ech cohomology
the element $1/x_1 \cdots x_d \in H^d_\fom (R)$,
and we call such a cohomology class a {\em canonical
element} or a {\em paraclass}.
If $R$ contains a field $K$, then such a paraclass is $\neq 0$
(The proof for this is different for positive and zero charcateristic,
see \cite[Remark 9.2.4]{brunsherzogrev}).
This is not known for mixed characteristic,
but at least there exist paraclasses $\neq 0$.

The main ingredient of the theory of solid closure is
the concept of a solid algebra.
For a domain $R$, an $R$-algebra $A$ is called {\em solid} if there exists
a non-trivial $R$-module homomorphism $\varphi: A \ra R$.
If $R$ is a local complete Noetherian domain, this is equivalent with
$H^d_\fom (A) \neq 0$, see \cite[Corollary 2.4]{hochstersolid}.
In the following definition we require the stronger
condition that all the paraclasses
coming from $H^d_\fom(R)$ shall not vanish in $H^d_\fom (A)$.

\begin{definition}
Let $R$ denote a local Noetherian ring of dimension $d$
and let $A$ denote an
$R$-algebra. $A$ is called {\em parasolid} if
the image of every paraclass $0 \neq c \in H^d_\fom(R)$ in
$H^d_\fom (A)$ does not vanish.
\end{definition}

\begin{definition}
Let $R$ denote a Noetherian commutative ring and let $A$ be an $R$-algebra.
We say that $A$ is {\em parasolid} if $A_\fom$ is parasolid
over $R_\fom$ for every maximal ideal $\fom \subset R$.

\smallskip
We say that $A$ is {\em universally parasolid}
(even if $R$ is not Noetherian)
if for every local Noetherian ring
$R \ra R'$ the algebra $A \otimes_R R'$ is a parasolid $R'$-algebra.
\end{definition}

\begin{remark}
The condition $c \neq 0$ is imposed only for mixed characteristic,
for this is not known in general.
The paraclass $c=1/x_1 \cdots x_d \in H^d_\fom (R)$
vanishes in $H^d_\fom(A)$ (we will often say that it vanishes in $A$)
if and only if
$(x_1 \cdots x_d)^k \in (x_1^{k+1}, \ldots, x_d^{k+1})$
holds in $A$ for some $k \in \NN$.

In the definition of universally parasolid it is enough to consider only
complete local rings, since completion does not change the $d$-th local
cohomology.
In equal characteristic we have only to look at complete normal
domains.
\end{remark}

\medskip
Some examples of parasolid algebras are gathered together in the following
proposition.

\begin{proposition}
\label{parasolidexam}
Let $R$ denote a commutative ring and let
$A$ denote an $R$-algebra.
In the following situations $A$ is universally parasolid.

\renewcommand{\labelenumi}{(\roman{enumi})}
\begin{enumerate}

\item
$R$ is a direct summand of $A$.

\item
$A$ is faithfully flat over $R$.

\item
$R$ contains a field and $R \subseteq A$ is a finite extension
{\rm (}or $R \ra A$ is finite and $\Spec\, A \ra \Spec\, R$
is surjective{\rm )}.

\item
There exists $A \ra A'$ such that $A'$ is a universally
parasolid $R$-algebra.

\end{enumerate}
\end{proposition}

\proof
We may assume that $R$ is a local Noetherian ring
of dimension $d$ and we have to show that $A$ is parasolid.

(i). Let $A=R \oplus V$. Then the cohomological mapping
$H^d_\fom(R) \ra H^d_\fom (A)= H^d_\fom (R) \oplus H^d_\fom (V)$ is injective.

(ii).
In the flat case we know that
$H^d_\fom (A)=  H^d_\fom (R) \otimes_R A $.
Since $A$ is faithfully flat, it is also pure and again the
mapping is injective.

(iii).
Let $x_1, \ldots ,x_d$ be parameters of $R$ and let
$(x_1, \ldots, x_d)A \subseteq \fon $ be a maximal ideal of $A$.
Then $x_1, \ldots ,x_d$ are parameters in $A_\fon$ and their
paraclass does not vanish.

(iv) is clear.
\qed

\begin{proposition}
\label{parasolidcompl}
Let $R$ denote a local Noetherian ring of dimension $d$
and let $\hat{R}$ be its completion.
Then an $R$-algebra $A$ is parasolid if and only if
$A'= A \otimes_R \hat{R}$ is parasolid over $\hat{R}$.
\end{proposition}

\proof
(Here I have to thank Prof. Storch for a hint)
We know that $H^d_\fom(R)= H^d_\fom (\hat{R})$ and hence the
mapping on the cohomology
$H^d_\fom (A) \ra H^d_\fom (A') = H^d _\fom (A) \otimes _R R'$
is injective. Hence it is enough to show that the paraclasses
coming from $R$ and from $\hat{R}$ are the same up to a unit.
Let $x_1, \ldots , x_d$ be parameters in $\hat{R}$ and let $\foa$
be the ideal they generate.
Since $\foa$ is primary, it is the extended ideal of an ideal $\fob$ in $R$.
The number of minimal generators does not change, hence $\fob$
is also a parameter ideal. Therefore $(x_1, \ldots ,x_d) = (y_1, \ldots ,y_d)$,
where $y_i \in R$.

Then they define up to a unit the same paraclass. This can be proved
by replacing inductively $x_i$ by $y_i$ using
\cite[Exercises 5.1.10-5.1.15]{brodmann}.
\qed

\medskip
Over a complete local domain a parasolid algebra is solid,
and in positive characteristic these notions coincide.

\begin{lemma}
\label{parasolidpos}
Let $K$ denote a field of positive characteristic $p$ and let
$R$ denote a local Noetherian $K$-Algebra of dimension
$d$. Let $A$ be an $R$-algebra.
Then the following are equivalent.

\renewcommand{\labelenumi}{(\roman{enumi})}
\begin{enumerate}

\item
$A$ is a parasolid $R$-algebra.

\item
There exists a paraclass $c \in H^d_\fom (R)$
which does not vanish in $H^d_\fom (A)$.

\item
$H^d_\fom(A) \neq 0$.

\end{enumerate}
\end{lemma}

\proof
(i) $\Rightarrow$ (ii) and (ii) $\Rightarrow $ (iii) are clear.
(iii) $\Rightarrow $ (i).
Suppose there exists a paraclass $1/(x_1 \cdots x_d)$ which does
vanish in $H^d_\fom (A)$,
and let $c= a/(x_1^r \cdots x_d^r)$ be an element in
$H^d_\fom(A)$.
The vanishing of the paraclass means that there exists an equation
in $A_{x_1 \cdots x_d}$,
$$ \frac{1}{x_1 \cdots x_d} \, = \, \frac{a_1}{(x_2 \cdots x_{d})^k} 
+ \ldots + \frac{a_d }{(x_1 \cdots x_{d-1})^k} \, .$$
Let $q=p^{e} \geq r$ and apply the $e$-th Frobenius to this equation.
This shows that also the paraclass $1/(x_1 \cdots x_d)^q $ vanishes,
hence also $c=0$.
\qed

\begin{corollary}
\label{parasolidpos2}
Let $K$ denote a field of positive characteristic $p$ and let
$R$ denote a local Noetherian $K$-Algebra of dimension
$d$. Let $A$ be an $R$-algebra.
Then the following are equivalent.

\renewcommand{\labelenumi}{(\roman{enumi})}
\begin{enumerate}

\item
$A$ is universally parasolid.

\item
$A \otimes_R \hat{R}/\foq$ is solid for every minimal prime of
the completion $\hat{R}$ {\rm(}his property is called formally solid,
see \cite[Definition 3.2]{hochstersolid}{\rm )}.

\end{enumerate}

If $R$ is a complete local domain, then $A$ is parasolid if and only
if it is universally parasolid.
\end{corollary}
\proof
(i) $\Rightarrow $ (i) is clear, so assume that (ii) holds.
Then from \cite[Theorem 3.7]{hochstersolid} we know that
$A'=A \otimes _R R'$ is solid for every Noetherian complete local domain $R'$
and \ref{parasolidpos} gives that $A'$ is parasolid.
\qed

\begin{example}
Let $R=K[[X,Y,Z]]/(XY,XZ)$, so that its spectrum
is the union of an affine plane and an affine line meeting
in one point. Let $A=K[[X,Y,Z]]/(X)$.
Then $A$ is a parasolid $R$-algebra,
but not universally parasolid. 
\end{example}

\medskip
The following proposition 
characterizes parasolid $R$-algebras $A$ for low dimensions ($\leq 2$)
by $H^d_\fom (A) \neq 0$.
Note that the characterization is set-theoretical for $d=0$,
topological for $d=1$ and algebro-geometrical for $d=2$.

\begin{proposition}
\label{paratwodimensional}
Let $R$ denote a local Noetherian ring of dimension $d$ with
maximal ideal $\fom$
and let
$A$ denote an $R$-algebra. Then the following are
equivalent.

$d=0$.
\renewcommand{\labelenumi}{(\roman{enumi})}
\begin{enumerate}

\item
$A$ is parasolid.

\item
$A \neq 0$.

\item
$\Spec\, A \neq \emptyset $.

\end{enumerate}

$d=1$. Let $x \in R$ be such that $V(x) = V(\fom)$.

\renewcommand{\labelenumi}{(\roman{enumi})}
\begin{enumerate}

\item
$A$ is parasolid.

\item
$H^1_\fom(A) \neq 0$.

\item
$A \lra A_x$ is not surjective
{\rm (}this means that $x^k \not \in (x^{k+1})$ in $A$
for all $k \in \NN$.{\rm )}

\item
$D(x)=D(\fom A) \subseteq \Spec\, A$ is not a closed subset.
\end{enumerate}

$d=2$.
\renewcommand{\labelenumi}{(\roman{enumi})}
\begin{enumerate}

\item
$A$ is parasolid.

\item
$H^2_\fom(A) \neq 0$.

\item
$D(\fom A) \subseteq \Spec\, A$ is not affine.
\end{enumerate}
\end{proposition}

\proof
Suppose $d=0$.
For every $R$-module $M$ we have
$$H^0_\fom (M) = \Gamma_\fom (M)= \{ v \in M:\, \exists \, n  \in \NN
\, \mbox{ such that }\, 
\fom^n v =0 \} =M \, .$$
In particular $H^0_\fom (R) = R$ and $1 \in R$ is the only paraclass,
and it vanishes in $H^0_\fom (A)$ if and only if $A=0$.
(ii) $\Leftrightarrow $ (iii) is clear.

Suppose $d=1$.
(i) $\Rightarrow $ (ii) is clear. The cohomology sequence is
$$0 \lra H^0_\fom (A) \lra A \lra A_x \lra H^1_\fom (A) \lra 0 \, $$
and this gives (ii) $\Leftrightarrow $ (iii).
Since this is true for every parameter $x \in R$,
we conclude (ii) $\Rightarrow $ (i).
If $A \ra A_x$ is surjective, then we have a closed embedding
$D(x) \ra \Spec\, A$.
On the other hand, consider
$A \ra A/\foa \ra A_x$, where $\foa$ is the kernel.
If $D(x)$ is closed, then
$D(x)=V(\foa)$, but $x$ is a unit in $A/\foa$,
hence $A/\foa  \cong A_x$.

Suppose $d=2$.
(i) $\Rightarrow $ (ii) is again clear.
If (ii) holds, then $H^1(D(\fom A), \O_A) \cong  H^2_\fom (A) \neq 0$
and $D(\fom A)$ is not affine.

(iii) $\Rightarrow$ (i).
Let $x,y$ be parameters in $R$ and assume that
$c=1/xy$ vanishes in $H^2_\fom (A) = H^2_{\fom A} (A)$.
This means that there exists an equation
$$ \frac{1}{xy} = \frac{a}{x^k} + \frac{b}{y^k} \, , \, \, \, a,b \in A .$$
Hence we see that $\, q = b/y^{k-1} = 1/x - ay/x^k \,$ and
$ \, p= a/x^{k-1} = 1/y - bx/y^{k} \,$ are functions defined on
$D(x,y)=D (\fom A) \subseteq  \Spec\, A$.
Furthermore we see that
$$ 1 = \frac{a}{x^{k-1}} y + \frac{b}{y^{k-1}}x =px+qy \, ,$$
and this shows that $D(x,y)$
is affine due to \cite[Exercise II.2.17]{haralg}.
\qed

\section{Order ideals}
\label{ordersec}
Let $A$ denote an $R$-algebra. An $R$-linear homomorphism
$\varphi \in \Hom_R(A,R)$ induce a mapping on the cohomology
$H^d_\fom (\varphi): H^d_\fom (A) \ra H^d_\fom (R)$.
These mappings are useful in
the study of $H^d_\fom (A)$, in particular if $R$ is complete.

\begin{definition}
Let $A$ denote an $R$-algebra. 
The ideal
$$U(A)
= \{ \varphi(1):\, \, \varphi: A \longrightarrow R\, \mbox{ $R$-linear }\}\, ,$$
is called the {\em order ideal} of the algebra $A$.
\end{definition}

\begin{remark}
It is clear that this is an ideal and that
an element $u \in R$ belongs to $U(A)$
if and only if $u$ belongs to the image of
some $R$-module homomorphism $A \ra R$.
The algebra $A$ contains $R$ as a direct summand if and only if
$1 \in U(A)$. The algebra $A$ over the domain $R$
is solid if and only if $U(A) \neq 0$.
\end{remark}

Let $R$ be a local Noetherian ring with residue class field
$k$ and let $E$ be an injective hull of $k$.
Consider the Matlis-functor
$D(-)=Hom_R(- ,E)$.
Let $M$ denote an $R$-module. Then the mapping
$$ H^d_\fom(M) \longrightarrow D(D(H^d_\fom(M)) $$
is injective due to \cite[Remarks 10.2.2]{brodmann}.
This means that for an element
$0 \neq c \in H^d_\fom(M)$
there exists an $R$-morphism
$ \psi: H^d_\fom(M) \longrightarrow E$ such that
$\psi(c) \neq 0$.

If $R$ is complete, then every $R$-module homomorphism
$E \longrightarrow E$ is the multiplication 
by a scalar $r \in R$, see \cite[Theorem 10.2.12]{brodmann}.
If furthermore $R$ is Gorenstein,
then $E=H^d_\fom(R)$ is an injective hull of the residue class field,
see \cite[Lemma 11.2.3]{brodmann}.

\begin{lemma}
\label{bijektivkoho}
Let $R$ denote a local complete Gorenstein ring of dimension $d$
and let $A$ be an $R$-algebra.
Then the mapping
$$\Hom_R (A,R) \lra \, Hom_R(H^d_\fom(A), H^d_\fom(R)),\, \, \, \, 
\varphi \longmapsto H^d_\fom(\varphi) $$
is bijective.
\end{lemma}

\proof
We define an inverse mapping.
Let $\psi :\, H_\fom^d(A) \ra H_\fom^d(R)$ be $R$-linear
and let $a \in A$.
Consider
$$H^d_\fom(R) \stackrel{i}{\lra}
H^d_\fom(A) \stackrel{a}{\lra} H^d_\fom(A)
\stackrel{\psi}{\lra} H^d_\fom(R) \, .$$
The whole mapping is multiplication by an element $r \in R$, and we
set $\varphi_\psi (a):=r$.
This gives an $R$-linear mapping $\varphi_\psi:\, A \ra R$.

Let $\varphi: A \ra R$ be given.
Then
$$(H_\fom^d(\varphi) \circ a \circ i) ( \frac{s}{x_1 \cdots x_d})
= H^d_\fom(\varphi) (\frac{as}{x_1 \cdots x_d})
=  \frac{ \varphi(a) s }{x_1 \cdots x_d} \, ,$$
so this is multiplication
by $\varphi (a)$.

Let $\psi: H^d_\fom(A) \ra H^d_\fom(R)$ be given and let
$\varphi _\psi$ as just defined.
Then
$$H^d_\fom (\varphi_\psi) (\frac{a}{x_1 \cdots x_d})
= \frac{\varphi_\psi(a)}{x_1 \cdots x_d}
=\psi ( \frac{a}{x_1 \cdots x_d}) \, .$$
\qed

\begin{lemma}
\label{classnonvanish}
Let $R$ denote a complete local Gorenstein ring of dimension $d$
and let $A$ be an $R$-algebra.
Let $c \in H^d_\fom(R)$. Then the following are equivalent.

\renewcommand{\labelenumi}{(\roman{enumi})}
\begin{enumerate}

\item 
$i(c) \neq 0$ in $H^d_\fom (A)$.

\item
There exists a homomorphism $\varphi: A \ra R$ such that
$0 \neq H^d_\fom(\varphi)(i(c)) \in H^d_\fom(R)$.

\item
There exists an element $u \in U(A)$ such that
$uc \neq 0$ in $H^d_\fom (R)$.

\item
$U(A) \not\subseteq \Ann_R \, c$.
\end{enumerate}
\end{lemma}

\proof
(i) $\Rightarrow $ (ii). From Matlis duality we know that there exists
$\psi: H^d_\fom (A) \ra H^d_\fom (R)$ such that $\psi (c) \neq 0$.
Due to \ref{bijektivkoho} this mapping is induced by
$\varphi: A \ra R$. 
Suppose (ii) holds. Then $u= \varphi(1)$ satisfies (iii).
From (iii) to (i) is clear, and (iii) and (iv) are the same.
\qed

\begin{corollary}
\label{parasolidorder}
Let $R$ denote a complete local Gorenstein ring of dimension $d$
and let $A$ be an $R$-algebra.
Then the following are equivalent.

\renewcommand{\labelenumi}{(\roman{enumi})}
\begin{enumerate}

\item 
$A$ is parasolid.

\item
$U(A)$ is contained in no annihilator of a paraclass.
\end{enumerate}
\end{corollary}

\proof
Follows from \ref{classnonvanish}.
\qed

\medskip
There exists always a class $c \in H^d_\fom (R)$ (in general no paraclass)
such that its annihilator is $\Ann_R \, c = \fom $.

\begin{corollary}
\label{directtest}
Let $R$ denote a complete local Gorenstein ring
and let $c \in H^d_\fom (R)$ be a class such that
$\Ann \, c = \fom $.
Let $A$ denote an $R$-algebra.
Then the following are equivalent.

\renewcommand{\labelenumi}{(\roman{enumi})}
\begin{enumerate}

\item
$R$ is a direct summand in $A$

\item
The mapping $i: H^d_\fom (R) \ra H^d_\fom (A)$ is injective.

\item
$i(c) \neq 0$ in $H^d_\fom (A)$.

\item
$U(A)=R$.

\end{enumerate}
\end{corollary}

\proof
(i) $\Rightarrow $ (ii) and (ii) $\Rightarrow$ (iii) are clear.
(iii) $\Rightarrow $ (iv). Let $i(c) \neq 0$.
From \ref{classnonvanish} we know that $U(A) \not\subseteq \Ann \, c = \fom$,
hence $U(A)$ contains a unit and then $U(A)=R$.
(iv) $\Rightarrow$ (i) is clear.
\qed

\begin{lemma}
\label{charactorder}
Let $R$ denote a local complete Gorenstein ring of dimension $d$
and let $A$ denote an $R$-algebra.
Let $i: H^d_\fom(R) \ra H^d_\fom(A)$ denote the mapping on the
local cohomology. Let $u \in R$.
Then the following are equivalent.

\renewcommand{\labelenumi}{(\roman{enumi})}
\begin{enumerate}

\item
$u \in U(A)$.

\item
$ u \in  \bigcap_{c \in H^d_\fom (R)}( \Ann_R \, c : \Ann_R \, i(c) ) \, .$

\item
$u$ annihilates $\ker \, i$.

\end{enumerate}
\end{lemma}

\proof
Suppose that (i) holds, say $u = \varphi(1)$.
Let $r \in R$ and suppose that
$r i(c)=0$ in $H^d_\fom (A)$.
Then $0=\varphi(r i(c))= u rc=$ and $ur$ annihilates $c$.

Suppose that (ii) holds, and let $i(c)=0$.
Then $\Ann \, i(c)= R$ and $\Ann \, (c) : R = \Ann \, (c)$, hence
$u \in \Ann \, (c)$.

Suppose that (iii) holds.
We have an injection $H^d_\fom(R)/ \ker\, i \hookrightarrow H^d_\fom (A)$.
Due to the condition multiplication on $H^d_\fom (R)$
by $u$ factors through $H^d_\fom(R)/ \ker\, i \ra H^d_\fom(R) $.
Since $H^d_\fom (R)$ is an injective module for a Gorenstein ring,
there exists also an extension
$H^d_\fom (A) \ra H^d_\fom (R)$.
This mapping comes from a mapping $\varphi: \,A \ra R$
and it must hold $\varphi(1)=u$.
\qed

\begin{lemma}
Let $R$ be a Noetherian local Cohen-Macaulay ring of dimension $d$
and let $x_1, \ldots , x_d$ be parameters.
Then the annihilator of the paraclass $1/x_1 \cdots x_d$ is
just $(x_1, \ldots , x_d)$.
\end{lemma}
\proof
Every system of parameters is
a regular sequence (\cite[Theorem 2.1.2]{brunsherzogrev}).
Consider the exact sequence
$0 \ra R \stackrel{x_d}{\ra} R \ra R/(x_d) \ra 0$.
The last part of the long exact sequence is
$$ 0 \lra H^{d-1}_\fom (R/(x_d)) \lra H^d_\fom (R)
\stackrel{x_d}{\lra} H^d_\fom (R) \lra 0 \, .$$
Here $1/x_1 \cdots x_{d-1} \mapsto 1/x_1 \cdots x_d$, and the
$R$-annihilator of these elements is the same.
Hence we do induction, where the beginning is clear, since
$H_\fom^0 (R) = R$ for $R$ zero-dimensional and
$0 =\Ann R$.
\qed

\medskip
The property that $R$ is a direct summand of $A$
is in general much stronger than the property that
$A$ is parasolid. For
complete regular rings however these properties coincide.
From this result we will deduce in section 4 that every ideal in a regular ring
is parasolidly closed.

\begin{theorem}
\label{regularparasolid}
Let $R$ be a complete regular local ring and let $A$ be an
$R$-algebra.
Then $A$ is parasolid if and only if $R$ is a direct summand of $A$.
\end{theorem}

\proof
Let $\fom=(x_1, \ldots, x_d)$, where $d$ is the dimension of $R$.
Then the annihilator of the corresponding paraclass is exactly
the maximal ideal $\fom$.
If $A$ is parasolid, then due to \ref{parasolidorder} we know that
$U(A) \not\subseteq \fom$, hence $U(A)=R$.
\qed

\section{parasolid closure}
\label{parasolidclosure}

Let $R$ denote a commutative ring and let $f_1, \ldots ,f_n$ and $h$
be elements in $R$. Then we call the $R$-algebra
$$A= R[T_1, \ldots, T_n]/(f_1T_1+ \ldots +f_nT_n+h)$$
the {\em forcing algebra} for the elements (or datas)
$f_1, \ldots, f_n;h$.
The forcing algebra $A$ forces $h \in (f_1, \ldots ,f_n)A$ and every
$R$-algebra $R \ra B$ such that $h \in (f_1, \ldots ,f_n)B$ factors
through $A$.

\begin{definition}
Let $R$ denote a commutative ring, let $f_1, \ldots ,f_n \in R$ and $h \in R$
be elements.
Let $A$ denote the forcing algebra for these elements.

If $R$ is local and Noetherian, then we say
that $h$ {\em belongs parasolidly} to $f_1, \ldots ,f_n$
if the forcing algebra $A$ is parasolid.

We say that $h$ {\em belongs universally parasolidly} to
$ (f_1, \ldots ,f_n) $ if
the forcing al\-ge\-bra $A$ is universally parasolid.
\end{definition}

\begin{lemma}
\label{basic}
Let $R$ be a local Noetherian ring of dimension $d$ and
let $f_1, \ldots, f_n$ and $h$ be elements.
Then the following hold.

\renewcommand{\labelenumi}{(\roman{enumi})}
\begin{enumerate}
\item 
If $h \in (f_1, \ldots , f_n)$, then it belongs
parasolidly to $f_1, \ldots ,f_n$.

\item
If $(f_1, \ldots ,f_n) \subseteq (g_1, \ldots ,g_m)$ and if
$h$ belongs parasolidly to $f_1, \ldots ,f_n$, then it belongs
parasolidly to $g_1, \ldots ,g_m$.

\item
If $h$ belongs parasolidly to $f_1, \ldots ,f_n$, then also
$rh$ belongs parasolidly to $f_1, \ldots ,f_n$.

\item
Let $a \in (f_1, \ldots ,f_n)$. Then $h$ belongs parasolidly to
$f_1, \ldots ,f_n$ if and only if $a+h$ belongs parasolidly to them.
\end{enumerate}
The same is true for $R$ any Noetherian ring if we replace
parasolidly by universally parasolidly.
\end{lemma}

\proof
(i).
Let $h=r_1f_1+ \ldots +r_nf_n$, $r_i \in R$. This defines a section
$A \ra R$, thus $R$ is a direct summand in $A$,
hence $A$ is parasolid.
(ii)
Let $f_i= r_{i1}g_1 + \ldots + r_{im}g_m$ for $i=1, \ldots ,n$.
This defines a mapping
$$R[S_1, \ldots , S_m]/(g_1S_1+ \ldots + g_mS_m +h)
\lra R[T_1, \ldots ,T_n]/(f_1T_1 + \ldots + f_nT_n +h)  $$
by sending
$ S_j \mapsto \sum_i r_{ij} T_i$.
This is well defined since
$$\sum_{j=1}^m g_j S_j +h
= \sum_{j=1}^m g_j (\sum_{i=1}^n r_{ij} T_i) +h
= \sum_{i=1}^n (\sum_{j=1}^m  r_{ij}g_j )T_i +h \, .$$
Since the algebra on the right is parasolid, so is the
algebra on the left.

(iii).
Here we use the mapping
$$R[S_1, \ldots , S_n]/(f_1S_1+ \ldots + f_nS_n +rh)
\lra R[T_1, \ldots ,T_n]/(f_1T_1 + \ldots + f_nT_n +h)  $$
which sends $S_i \ra r T_i$.

(iv).
Let $a=a_1f_1+ \ldots +a_nf_n$. Then we use the mapping
$$R[S_1, \ldots , S_n]/(f_1S_1+ \ldots +f_nS_n +h)
\lra R[T_1, \ldots ,T_n]/(f_1T_1 + \ldots + f_nT_n + a+h)  $$
where $S_i \mapsto T_i+a_i $.

The corresponding statement for universally parasolidly follows.
\qed

\begin{definition}
\label{parasolidcldef}
Let $R$ denote a local complete Noetherian domain and let
$h \in R$ and $ f_1, \ldots ,f_n \in R$ be elements.
We say that $h$ belongs to the {\em parasolid closure} of
$f_1, \ldots ,f_n$ if there exists a chain of elements
$f_{n+1}, \ldots ,f_m =h$
such that $f_i $ belongs universally parasolidly to
$f_1, \ldots , f_n, \ldots ,f_{i-1} $ for $i=n+1, \ldots, m$.
\end{definition}

\begin{lemma}
Let $R$ denote a local complete Noetherian domain and
let $f_1, \ldots, f_n$ be elements.
The elements which belong to the parasolid closure of
$f_1, \ldots ,f_n$ form an ideal, which depends
only on the ideal $(f_1, \ldots ,f_n)$.
\end{lemma}

\proof
Let $h$ belong to the parasolid closure.
This means that there exists a chain
$f_1, \ldots ,f_n,\, f_{n+1}, \ldots ,f_m=h$, where
$f_i $ belongs universally parasolidly to
$f_1, \ldots, f_n, \ldots, f_{i-1}$
for $i=n+1, \ldots ,m$.
This is also true if we replace $h$ by $rh$ due to \ref{basic}(iii).
Let $g$ be another element which belongs also to the parasolid closure,
and let $f_1, \ldots, f_n, g_1, \ldots, g_k=g$ be a corresponding
chain.
Then also $f_1, \ldots,f_n, f_{n+1}, \ldots, f_m,g_1, \ldots ,g_k$ is
such a chain due to \ref{basic} (ii), and then due to \ref{basic} (iv)
also $f_1, \ldots ,f_n, f_{n+1}, \ldots, f_m,g_1, \ldots ,g_k+f_m=g+h$
is such a chain.
\qed

\begin{definition}
\label{parasolidclosuredef}
Let $R$ denote a Noetherian ring and let $I \subseteq R$
denote an ideal.

If $R$ is a local complete domain, then we call
the ideal just defined the {\em parasolid closure} $I^\pasoclo$ of $I$.

For general $R$, we define $I^\pasoclo $ by the condition that
$h \in I^\pasoclo$ if and only if this is true
for all $\hat{R_\fom}/\foq$, where $\fom$ is a maximal ideal $\fom$ of $R$
and $\foq$ is a minimal prime of the completion
$\hat{R_\fom}$
(These rings are called the local complete domains of $R$).
\end{definition}

\begin{proposition}
\label{properties}
Let $R$ denote a Noetherian ring and let
$I \subseteq R$ be an ideal, $I=(f_1, \ldots, f_n)$.
Let $h \in R$.
Then the following hold.

\renewcommand{\labelenumi}{(\roman{enumi})}
\begin{enumerate}

\item
If $h $ belongs universally parasolidly to $I$, then $h \in I^\pasoclo$.

\item
$I \subseteq I^\pasoclo $

\item
If $I \subseteq J$, then $I^\pasoclo \subseteq J^\pasoclo$.

\item
The parasolid closure is persistent:
If  $h \in I^\pasoclo$, then for every
ring homomorphism $R \ra S$, $S$ Noetherian, it holds that
$h \in (IS)^\pasoclo $.

\item
Suppose that $R$ is a local complete domain.
Then $I=I^\pasoclo$ if and only if
no $h \not \in I$
belongs universally parasolidly to $(f_1, \ldots ,f_n)$.

\item
$(I^\pasoclo)^\pasoclo = I^\pasoclo $

\item
$h \in I^\pasoclo$ if and only if $h \in (IR_\fom)^\pasoclo$
for every localization at a maximal ideal {\rm (}or at a prime ideal{\rm )}.

\item
$h \in I^\pasoclo$ if and only if $h \in (IR/\foq)^\pasoclo$
for every minimal prime ideal $\foq$ of $R$.

\item
Suppose that $R$ is local. Then $h \in I^\pasoclo$ if and only if
$h \in (I\hat{R})^\pasoclo$.

\end{enumerate}
\end{proposition}

\proof
(i) and (ii) are clear.
In (iii) we may assume that $R$ is a local complete domain
and then thid follows from \ref{basic}(ii).
(iv). Let $S \ra S'$ be a complete local domain of $S$.
This mapping factors then through a complete local domain $R'$ of $R$,
and in $R'$ there exists a chain for $h$ as in the definition
\ref{parasolidcldef} and this gives the chain in $S'$.
(v) is clear, since every chain must start with an element
which belongs universally parasolidly to the ideal.
(vi). Let $h \in (I^\pasoclo)^\pasoclo$ and let
$R'$ denote a local complete domain of $R$.
Then $(I^\pasoclo)R' \subseteq (IR')^\pasoclo$ due to persistence
and then
$h \in (I^\pasoclo R')^\pasoclo \subseteq ( (IR')^\pasoclo )^\pasoclo$.
Hence we may assume that $R$ is a local complete domain,
and the statement is clear from the definition by chains.
(vii),(viii) and (ix) are clear by definition.
\qed

\begin{remark}
\label{variants}
A natural question here is: can we through away some of the burden
in the definition of parasolid closure which we needed to obtain
the basic properties listed in \ref{properties}?
This refers in particular to
the chains in \ref{parasolidcldef},
to localization and completion in \ref{parasolidclosuredef}
and to the use of all local Noetherian rings in the definition
of universally parasolid.
As for the last point, we will show
in section \ref{expresssec} that for $K$-algebras of finite type
it is enough to consider only local algebras which are
essentially of finite type.
(A variant with some advantages is
to consider only those local rings where
the extended ideal of $I$ is primary to the maximal ideal.)

\end{remark}

\begin{lemma}
\label{contraction}
Let $K$ be a field and let $R$ be a Noetherian $K$-algebra,
let $f_1, \ldots ,f_n \in R$ and $h \in R$.
Let $\varphi : R \ra S$ be a finite extension such that 
$h \in (f_1, \ldots ,f_n)S$.
Then $h \in (f_1,\ldots, f_n)^\pasoclo$.
\end{lemma}

\proof
Due to \ref{parasolidexam} (iii) the finite algebra
$S$ is universally parasolid.
\qed

\begin{remark}
The statement in the last lemma is open in the case of mixed characteristic
and equivalent to the monomial conjecture.
If we would define a parasolid algebra by the condition that
not the paraclasses $c \neq 0$ do not vanish,
but only these paraclasses which stay $\neq 0$ in every finite extension,
then the contraction property would also hold in mixed characteristic,
but then we could not prove the regular property anymore.
\end{remark}

\begin{example}
\label{nilpotent}
The parasolid closure of $0$ consists of the nilpotent elements.
One inclusion is clear from the definition which refers to domains.
If $h$ is not nilpotent, then there exists
a $R \ra R' \neq 0$ where $h$ becomes a unit
and then the forcing algebra $R'/(h)=0$ is not parasolid.

\end{example}

\begin{example}
Let $C$ denote a local three-dimensional regular ring with
maximal ideal $\fom =(x,y,z)$.
Let
$$R=C/(x^{i}+y^{j} -z^k), \,\  i,j,k \geq 1 \, .$$
$x$ and $y$ are then parameters in $R$,
we consider the forcing algebra for $x$, $y$ and $z$, hence
$$A=R[T_1,T_2]/(xT_1+ yT_2-z) \, .$$
We have then the equation
$$x^{i} +y^{j} =(xT_1+yT_2)^k \,\,\, \,  (\ast) \, \,  .$$

\smallskip
Suppose that $k > i$ (or $k >j$).
From the last equation we can read off the containment
$x^{i}(1-x^{k-i}T_1^k) \in (y)$.
Since $k >i$ we have $x^k \in (x^{i+1})$, hence we get
$x^{i} \in (y,x^{i+1})$ and this shows that $D(x,y) \subset \Spec\, A$
is affine.
Then the forcing algebra is not parasolid (not even solid) and
$z \not\in (x,y)^\pasoclo$.

\smallskip
Suppose that $k =i =j$. From the equation $(\ast)$ we see
$x^k(1-T_1^k) \in (y)$, $y^k(1-T_2^k) \in (x)$ and
$ \left(
\begin{array}{c} 
k \cr 
r
\end{array}
\right)
x^ry^{k-r}T_1^r T_2^ {k-r} \in (x^{r+1},y^{k-r+1})$ for $1 \leq r \leq k-1$.
Suppose that $k \geq 2$ and that at least one of the
coefficients
$ \left(
\begin{array}{c} 
k \cr 
r
\end{array}
\right)$,
$1 \leq r \leq k-1$ is a unit in $R$.
Then again $x^ny^m \in (x^{n+1}, y^{m+1})$ for some $n, m$ and
$z \not\in (x,y)^\pasoclo$.
\end{example}

\section{Regular rings}
\label{regularsec}

\begin{theorem}
\label{regular}
Let $R$ be a regular local ring and let $I \subseteq R$ be an ideal.
Then $I=I^\pasoclo$.
\end{theorem}

\proof
We may assume that $R$ is complete.
Let $I=(f_1, \ldots ,f_n)$ and let $h \in R$ and let
$A$ denote the forcing algebra.
Suppose that $h$ belongs universally parasolidly to $f_1, \ldots ,f_n $.
This means in particular that $A$ is parasolid.
Hence due to \ref{regularparasolid} $R$ is a direct summand of $A$
and then $h \in (f_1, \ldots, f_n)$.
\qed

\medskip
We give a second (more constructive)
proof for the case that $R$ contains a field.
Again we may assume that $R$ is complete, hence
$R =K[[X_1, \ldots ,X_d]]$.
Assume that $h \not\in (f_1, \ldots, f_n)$.
Then also
$h \not\in (f_1, \ldots ,f_n ,X_1^m, \ldots , X_d^m)$
for $m \in \ZZ$ big enough and then we find
$r=(r_1,\ldots ,r_d)$ such that
$h \not\in (I, X_1^{r_1+1},\ldots ,X_d^{r_d+1})$, but
$ h \in (I, X_1^{r_1},X_2^{r_2+1},\ldots ,X_d^{r_d+1})$, $\ldots $,
$h \in (I, X_1^{r_1+1},\ldots ,X_{d-1}^{r_{d-1}+1},X_d^{r_d}) $.
In particular we have
$hX_i \in (I, X_1^{r_1+1},\ldots ,X_d^{r_d+1})$ for
every variable $X_i$.
We may replace $I$ by
$(I, X_1^{r_1+1}, \ldots , X_d^{r_d+1})
=(f_1, \ldots, f_n,X_1^{r_1+1}, \ldots , X_d^{r_d+1}) $
and we may assume that
$f_j =  \sum_{\nu \leq r} f_{j,\nu} X^\nu$ and also
$h= \sum_{\nu \leq r} h_\nu X^\nu$.

There exists a $K$-linear form $\varphi :K[[X_1,\ldots ,X_d]] \ra K$
given by $X^\sigma \mapsto c_\sigma$ such that
$\varphi (I)=0$, but $\varphi(h)=1$.
Modulo $(X_i^{r_i+1})$ we have
$X^\sigma f_j= \sum_{\mu \leq r} f_{j,\mu - \sigma } X^\mu \in I$,
hence
$ \sum_{\mu \leq r} c_\mu f_{j,\mu - \sigma } = 0$
(where $f_{j, \lambda}=0$ if some entry of $\lambda $ is negative).
Furthermore
$X^\sigma h = \sum_{\mu \leq r} h_{\mu -\sigma} X^\mu $ belongs also to
$I$ for $\sigma >0$, hence
$ \sum _{\mu \leq r} c_\mu h_{\mu -\sigma} =0 $
for $ \sigma >0$, while $\sum_\mu c_\mu h_\mu =1$.
We consider
$$\sum_{\sigma \leq r} c_{r - \sigma} X^\sigma (f_1T_1+\ldots +f_nT_n+h) 
= \sum_{\nu}a_\nu X^\nu \, ,$$
where $a_\nu \in K[T_1,\ldots ,T_n]$ is the coefficient for $X^\nu$.
Then for $\nu \leq r$ we may write
\begin{eqnarray*}
a_\nu &= &
\sum_{\sigma} c_{r-\sigma}(f_1T_1+\ldots +f_nT_n+h)_{\nu - \sigma} \\
&=& \sum_\mu c_\mu (f_1T_1+\ldots +f_nT_n+h)_{\mu-(r-\nu)} \\
&=& (\sum_\mu c_\mu f_{1,\mu - (r-  \nu)})T_1
+\ldots +(\sum_\mu c_\mu f_{n,\mu - (r-  \nu)})T_n
+ \sum_\mu c_\mu h_{\mu - (r-  \nu)} \, \, .
\end{eqnarray*}
Since
$\sum_{\mu \leq r} c_\mu f_{j,\mu-(r- \nu)} =0$
for $\nu \leq r$ and
$\sum_{\mu} c_{\mu} h_{\mu-(r-\nu)}=0$ for $\nu < r$
and $\sum_\mu c_\mu h_\mu=1 $ for $\nu=r$ we see that
$a_\nu=0 $ for $\nu < r$ and $a_r=1$. Hence we get
$$ \sum_{\sigma} c_{r - \sigma} X^\sigma
(f_1T_1+\ldots +f_nT_n+h)=X_1^{r_1} \cdots X_d^{r_d} +R$$
where $R \in (X_1^{r_1+1},\ldots ,X_d^{r_d+1})$, which shows that
the paraclass for $X_1, \ldots, X_d$ vanishes in the forcing algebra.
\qed

\medskip
As usual in a tight closure type theory, we make the following
definition.

\begin{definition}
We call a Noetherian ring $R$ {\em pararegular} if $I=I^\pasoclo$ holds
for every ideal $I \subseteq R$.
\end{definition}

\medskip
The Theorem tells us that regular local rings are pararegular.
This is also true for regular rings in general, since a
locally pararegular ring is pararegular (the convers is not clear). 
In connection with
persistence we get the
following two corollaries.

\begin{corollary}
\label{pararegularpure}
Let $R$ be a Noetherian ring which is a pure subring
{\rm(}e. g. a direct summand{\rm)} of a
{\rm(}para{\rm )}regular ring $S$.
Then $R$ is pararegular.
\end{corollary}

\proof
Let $I \subseteq R$ denote an ideal, and let $h \in I^\pasoclo$.
Then also $h \in (IS)^\pasoclo =IS$, but $R \cap IS = I$.
\qed

\begin{corollary}
Let $R$ be a Noetherian ring and let $I \subseteq R$
be an ideal.
Then the parasolid closure of $I$ is contained in the
regular closure and in the integral closure,
$I^\pasoclo \subseteq I^{\rm reg} \subseteq  \bar{I}$.
\end{corollary}
\proof
Let $h \in I^\pasoclo$ and let $R \ra S$ be a regular local ring.
Then $h \in (IS)^\pasoclo =IS$, so by definition it lies inside
the regular closure.
\qed

\medskip
Of course one hopes that pararegular rings are reduced, normal
and Cohen-Macaulay, but we can prove only partial results.
The Cohen-Macaulayness of parasolid rings
follows in equal characteristic from the comparision with tight closure,
see section \ref{generic}.

\begin{lemma}
\label{principal}
Let $R$ denote a Noetherian domain and let $\tilde{R}$
denote its normalization.
Let $0 \neq f \in R$.
If $h \in (f)^\pasoclo$, then $h/f \in \tilde{R}$.
If $R$ contains a field, then also the converse holds.
\end{lemma}

\proof
Let $h \in (f)^\pasoclo$.
This is then true in all discrete valuation domains
$B$ of $R$ and hence $h \in fB$, so $h/f \in B$.
Then $h/f$ lies in the intersection of
all valuation domains and hence in the normalization.

If $h/f \in \tilde{R}$, then consider
the inclusion $R \subseteq R[h/f] \subseteq \tilde{R}$
and we find $h \in (f)$ in the finite extension $R[h/f]$,
which is universally parasolid if $R$ contains a field.
\qed

\begin{proposition}
Let $R$ denote a Noetherian ring and suppose that
$R$ is pararegular. Then the following holds.

\renewcommand{\labelenumi}{(\roman{enumi})}
\begin{enumerate}

\item
$R$ is reduced.

\item
If $R$ is a domain and contains a field, then it is normal.

\end{enumerate}
\end{proposition}
\proof
This follows from \ref{nilpotent} and \ref{principal}.
\qed

\begin{example}
Let $R=K[U,V,W]/(UV-W^2)=K[X^2,Y^2,XY] \subset K[X,Y]=S$.
Then the cohomology class $c= W/UV $ maps to $1/XY$.
Therefore $c \neq 0$, and its annihilator is the maximal ideal $(U,V,W)$.
The class $c$ is no paraclass in $H^2_\fom (R)$, but in $H^2_\fom (S)$.
Let $A$ denote an $R$-algebra where $i(c)=0$. This means that we have an
equation $W (UV)^{k} =T_1U^{k+1} + T_2 V^{k+1}$
and $B$ factors through
the forcing algebra $A$ for this equation.
This gives on $S$ the forcing algebra with equation
$XY (X^2Y^2)^{k} = T_1X^{2(k+1)} + T_2Y^{2(k+1)}$.
We know that $H^2_\fom (A \otimes_R S)=0$ due to
\ref{regular} and \ref{paratwodimensional}
and hence this is also true for $H^2_\fom (A)$, since $A \subset A \otimes_RS$
is a direct summand.
\end{example}

\section{solid closure and parasolid closure}
\label{solidparasolidsec}

We recall the definition of solid closure. Let
$R$ denote a Noetherian ring and let
$(f_1, \ldots, f_n)=I \subseteq R$ be an ideal, $h \in R$.
Then $h \in I^\soclo$ if and only if the forcing algebra $A$ on
$\hat{R_\fom}/\foq$ is solid,
where $\fom $ runs through the maximal ideals of $R$ and
$\foq$ is a minimal prime in the completion
$\hat{R_\fom}$.

\begin{proposition}
Let $R$ be a Noetherian ring and let $I \subseteq R$ be an ideal.
Then $I^\pasoclo \subseteq I^\soclo$.
\end{proposition}

\proof
We may assume that $R$ is a local complete domain.
It is enough to show that an element $h$ which belongs universally
parasolidly to $I$ belongs to the solid closure of $I$.
The forcing algebra $A$ for $I,h$ is then universally parasolid
and so in particular solid.
\qed

\begin{proposition}
Let $K$ denote a field of positive characteristic $p$ and let
$R$ denote a Noetherian $K$-Algebra.
Let $I=(f_1, \ldots ,f_n) \subseteq R$
denote an ideal. Then $I^\pasoclo = I^\soclo $.
\end{proposition}

\proof
We have already shown the direct inclusion, thus suppose
that $h \in I^\soclo$. Since solid closure is persistent
\cite[Theorem 5.6]{hochstersolid},
we know for every local Noetherian complete domain $R \ra R'$ that
$H^d_\fom (A') \neq 0$, where $d$ is the dimension of $R'$
and $A'$ is the forcing algebra over $R'$.
By \ref{parasolidpos} $A'$ is parasolid.
\qed

\medskip
We obtain the connection to tight closure in positive characteristic.

\begin{corollary}
\label{tightpasoclop}
Let $K$ denote a field of positive characteristic $p$ and let
$R$ denote a Noetherian $K$-Algebra.
Suppose that $R$ is an algebra essentially of finite type
over an excellent local ring or that the Frobenius endomorphism in $R$
is finite.
Then $I^* = I^\pasoclo = I^\soclo $ holds for every ideal $I \subseteq R$.
\end{corollary}

\proof
The equation $I^* =I ^\soclo $ is due to \cite[Theorem 8.6]{hochstersolid}.
\qed

\begin{remark}
\label{tightpasoclop2}
The conditions imposed in the corollary imply that $R$ has
a completely stable weak test element
for tight closure. Even without this condition it is true
that the tight closure lies inside the solid closure (= parasolid closure).
\end{remark}

\begin{example}
We have a look at the example of Roberts.
Let $K$ be a field of characteristic zero and consider in
the polynomial ring
$R=K[X,Y,Z]$ the elements $X^3,Y^3,Z^3$ and $X^2Y^2Z^2$.
The forcing algebra is then
$$ A=K[X,Y,Z] [T_1, T_2,T_3] /(X^3T_1 + Y^3T_2 +Z^3T_3+X^2Y^2Z^2) \, .$$
This algebra is not parasolid, since the paraclass $1/XYZ$ does vanish
in $H^3_\fom (A)$, and $X^2Y^2Z^2$ does not belong to
$(X^3,Y^3,Z^3)^\pasoclo =(X^3,Y^3,Z^3)$.

On the other hand, the computation of Roberts in \cite{robertscomp} shows
that the class $1/X^2Y^2Z^2$ does not vanish in $H^3_\fom (A)$ and
so $A$ is a solid algebra (at least over $\hat{R}$) and
$X^2Y^2Z^2$ belongs to $(X^3,Y^3,Z^3)^\soclo$.
This example also shows that some paraclasses may survive, while others
do not.
\end{example}

\section{Order ideal and test ideal}
\label{testandorder}

\medskip
An $R$-module homomorphism
$A=R[T_1, \ldots ,T_n]/(f_1T_1+ \ldots +f_nT_n+h) \ra R$
is given by $T^\nu \mapsto r_\nu \in R$ such that for every $\nu$
the condition
$$f_1 r_{\nu+e_1} + \ldots +f_n r_{\nu +e_n} = r_\nu $$
holds, where $e_i$ is a standard base vector.
It is in general not easy to construct such a morphism.
The elements $r_0$ of such module homomorphisms build
the order ideal.
We can use such order elements to compute tight closure.

\begin{lemma}
Let $R$ be a commutative ring and let
$A$ be a forcing algebra for the elements $f_1, \ldots ,f_n,h$.
Let $u \in U(A)$.
Then $uh \in (f_1, \ldots , f_n)$.

If furthermore $R$ contains a field of positive characteristic $p$,
then $uh^q \in I^{[q]}$ for all $q=p^{e}$.
\end{lemma}

\proof
Let $u= \varphi(1)$. We apply $\varphi$ to the forcing equation
and get
$$ uh= \varphi(1) h=\varphi(h) 
= f_1 \varphi(T_1) +\ldots + f_n \varphi(T_n) \, .$$
Under the Frobenius homomorphism the forcing equation becomes
$f_1^qT_1^q+ \ldots +f_n^qT_n^q = h^q$. We apply again $\varphi$
and get $ uh^q = f_1^q \varphi (T_1^q) + \ldots + f_n^q \varphi (T_n^q)$.
\qed

\begin{remark}
If $R$ is a domain over a field of positive characteristic,
and if $U(A) \neq 0$, then
an element $u \neq 0$ can be taken to show that $h \in I^*$.
It is not clear whether an element $u$ such that
$uh^q \in I^{[q]}$ for all $q=p^{e}$ does belong
to the order ideal (in the complete case).
The proof of \cite[Theorem 8.6]{hochstersolid}
seems to give only that a power of $u$ is an order element.
\end{remark}

\medskip
The following ideal is a candidate for a test ideal for parasolid
closure.

\begin{definition}
Let $R$ denote a local Noetherian complete ring of dimension $d$.
We set
$$U(R):= \bigcap_{A\mbox{ universally parasolid forcing algebra}} U(A) \, .$$
and
$$V = \langle c \in H^d_\fom (R)\,: \exists \, \,
\mbox{universally parasolid forcing algebra } i:R \ra A, \,
i(c)=0  \rangle \, .$$

\end{definition}

\begin{lemma}
\label{uannv}
Let $R$ denote a local complete Gorenstein ring.
Then $U(R) =\Ann \, V$.
\end{lemma}
\proof
This follows from \ref{charactorder}. Let $u \in U(R)$ and let
$c \in V$ such that $c$ vanishes in the universally parasolid
forcing algebra $A$. By definition $u \in U(A)$ and then
$uc =0$, and this holds also in the submodule generated by such $c$.
On the other hand, suppose that $u \in \Ann \, V$ and let
$i: R \ra A$ be a universally parasolid forcing algebra.
Then $\ker \, i \subseteq V$ and $u$ annihilates the kernel,
hence $u \in U(A)$. (This correspondence is true for every family
of algebras $A_i$.)
\qed

\begin{proposition}
Let $R$ denote a local complete Gorenstein ring
which contains a field of positive characteristic.
Let $T$ be the tight closure test ideal of $R$. Then
$T=U(R)$.
\end{proposition}
\proof
We know that
$T= \Ann_R( 0^*_{H_\fom^d(R)} )$, where
$0^*_{H_\fom^d(R)}$ is the tight closure of the submodule
$0 \subseteq H^d_\fom (R)$, see \cite[Proposition 4.1]{hunekeparameter}.
Due to \ref{uannv} it is enough to show that $V= 0^*_{H^d_\fom (R)}$.
We know
$0^*_{H_\fom^d(R)}
=\{ s/(x_1 \cdots x_d)^k:\, s \in (x_1^k, \ldots, x_n^k)^* \}$.
The inclusion
$0^*_{H_\fom^d(R)} \subseteq V$ is clear, since 
$ s/(x_1 \cdots x_d)^k$ is zero in the forcing algebra
for $s;x_1^k, \ldots,x_d^k$ and this forcing algebra is
(formally solid, hence) universally parasolid.

So suppose that
$c = s/(x_1 \cdots x_d)^k$ vanishes in a
universally parasolid algebra $B$. This means that
there exists $m \geq 1$ such that
$s (x_1 \cdots x_d)^m \in (x_1^{k+m}, \ldots ,x_d^{k+m} )$ in $B$.
Then the forcing algebra $A$ for the elements
$s(x_1 \cdots x_d)^m; x_1^{k+m}, \ldots , x_d^{k+m} $ maps to $B$,
and since $B$ is universally parasolid also $A$ is universally parasolid.
Hence
$s(x_1 \cdots x_d)^m \in (x_1^{k+m}, \ldots ,x_d^{k+m} )^\pasoclo
= (x_1^{k+m}, \ldots ,x_d^{k+m} )^*$
and this means that
$s(x_1 \cdots x_d)^m /(x_1^{k+m} \cdots x_d^{k+m} ) = s/(x_1^k \cdots x_d^k)
\in 0^*$.
Therefore we have $0^*_{H_\fom^d(R)} = V$.
\qed

\section{Pararegular rings}
\label{pararegularsec}

For a complete local Gorenstein ring
there are a lot of different characterizations to be
pararegular.

\begin{proposition}
\label{pararegularchar}
Let $R$ denote a local complete Gorenstein ring of dimension $d$.
Let $c \in H^d _\fom(R)$ be a class such that $\Ann\, c= \fom$
{\rm(}only in part {\rm(}vii{\rm)}{\rm )}.
Then the following are equivalent.

\renewcommand{\labelenumi}{(\roman{enumi})}
\begin{enumerate}

\item
$R$ is pararegular, i.e. every ideal $I \subseteq R$ is parasolidly closed.

\item
Every universally parasolid forcing algebra $A$ contains $R$
as a direct summand.

\item
Every ideal generated by parameters is parasolidly closed.

\item
$V=0$.

\item
$U(R)=R$.

\item
Every universally parasolid $R$-algebra $B$ contains $R$ as a direct summand.

\item
$i(c) \neq 0$ in every universally parasolid
$R$-algebra $B$.

\end{enumerate}
\end{proposition}

\proof
(i) and (ii) are equivalent for every local complete Noetherian ring
due to \ref{properties} (v). (i) $\Rightarrow$ (iii) is clear.

(iii) and (iv) are equivalent for a local complete Cohen-Macaulay ring.
Suppose that (iii) holds. Let $x_1, \ldots ,x_d$ be parameters
and let
$c=s/x_1 \cdots x_d$ and suppose that $i(c)=0$ vanishes in the
universally parasolid forcing algebra $A$.
This means that
$s(x_1 \cdots x_d)^k \in (x_1^{k+1}, \ldots ,x_d^{k+1})$ holds in $A$
and then $A$ factors through the forcing algebra to
$x_1^{k+1}, \ldots ,x_d^{k+1}; s(x_1 \cdots x_d)^k$,
which is then also universally parasolid.
This means that
$s(x_1 \cdots x_d)^k \in (x_1^{k+1}, \ldots ,x_d^{k+1})^\pasoclo
= ( x_1^{k+1}, \ldots ,x_d^{k+1})$ in $R$ and hence $c=0$.
For the converse suppose that $s$ belongs universally
parasolidly to $(x_1, \ldots, x_d)$.
Then $c=s/x_1 \cdots x_d$ is zero in the universally parasolid
forcing Algebra, hence $c \in V=0$.
Then $s$ annihilates the paraclass $1/x_1 \cdots x_d$
and since $R$ is Cohen-Macaulay,
it follows that $s \in (x_1, \ldots ,x_d)$.

Suppose furthermore that $R$ is a complete Gorenstein ring.
The equivalence (iv) $\Leftrightarrow $ (v) follows from \ref{uannv}.

(iii) $\Rightarrow $ (vi).
Let $i: R \ra B$
be a universally parasolid algebra.
Due to \ref{charactorder}
it is enough to show that
the cohomological mapping
$H^d_\fom (R) \ra H^d_\fom(B)$ is injective.
So let
$c=s/x_1 \cdots x_d$ and suppose that $i(c)=0$.
This means that
$s(x_1 \cdots x_d)^k \in (x_1^{k+1}, \ldots ,x_d^{k+1})$ holds in $B$
and then again
$s(x_1 \cdots x_d)^k \in (x_1^{k+1}, \ldots ,x_d^{k+1})^\pasoclo
= ( x_1^{k+1}, \ldots ,x_d^{k+1})$ in $R$, thus $c=0$.

(vi) $\Rightarrow $ (vii) and (vi) $\Rightarrow $ (ii) are clear.
(vii) $\Rightarrow$ (v). If $A$ is a universally parasolid
forcing algebra, then $i(c) \neq 0$,
hence
$U(A) \not \subseteq \Ann\, c = \fom$, thus $U(A)=R$.
\qed

\begin{corollary}
Let $R$ denote a local Gorenstein ring.
Then $R$ is pararegular if and only if $\hat{R}$ is pararegular.
\end{corollary}

\proof
$ \Leftarrow $ is always true and follows from \ref{pararegularpure}.
The converse holds because of \ref{pararegularchar} (iii).
\qed

\begin{remark}
Since direct summands of regular rings are pararegular,
we get a lot of pararegular but non-regular rings.
If we work over a base field $K$, then all the quotient singularities
are pararegular, for example $K[x,y,z]/(xy-z^n)$.
This equation yields also in mixed characteristic
pararegular rings, as the following computation shows.
\end{remark}

\begin{example}
\label{toricexample}
Let $C$ denote a local three-dimensional regular ring
with maximal ideal $\fom =(x,y,z)$
and consider 
$$R=C/(xy-z^n), \,\,  n \geq 1 \, .$$
The elements $x$ and $y$ are now parameters in $R$.
We want to show that $R$ is pararegular.
The annihilator of
$c=z^{n-1}/xy \in H^2_\fom (R)$ is the maximal ideal $(x,y,z)$.
Due to \ref{pararegularchar}
we have to show that if $c$ vanishes in an $R$-algebra $A$,
then $A$ is not universally parasolid, and it is enough to look at
a forcing algebra $A$ given by an equation
$$z^{n-1} (xy)^k =T_1x^{k+1} +T_2y^{k+1}, \,  \,\, k \geq 0 .$$
This equation yields
$$1=\frac{T_1x^{k+1}}{z^{n-1} (xy)^k} + \frac{T_2y^{k+1}}{z^{n-1}(xy)^k}
=\frac{T_1x}{z^{n-1}y^k} +\frac{T_2y}{z^{n-1}x^k} \, \, \, \, \, (*)\, \,  .$$
The equations
$$ \frac{T_1z}{y^{k+1}} =\frac{T_1xz}{xy^{k+1}} 
= \frac{T_1xz}{z^{n}y^k} =\frac{T_1x}{z^{n-1}y^k}
= 1- \frac{T_2y}{z^{n-1}x^k}=1- \frac{T_2z}{x^{k+1}}$$
show that $T_1x/z^{n-1}y^k$ (and also $T_2y/z^{n-1}x^k$)
is a function defined on $D(x,y) \subseteq \Spec \, A$.
We take the $2n$ power of $(*)$ and get
$$ 1 =(\frac{T_1x}{z^{n-1}y^k} + \frac{T_2y}{z^{n-1}x^k})^{2n}
=P(\frac{T_1x}{z^{n-1}y^k})^n + Q(\frac{T_2y}{z^{n-1}x^k})^n \, ,$$
where $P$ and $Q$ are polynomials in $T_1x/z^{n-1}y^k$ and $T_2y/z^{n-1}x^k$
and hence defined on $D(x,y)$.
We write
$$  (\frac{T_1x}{z^{n-1}y^k})^n
= \frac{T_1^n x^n}{z^{(n-1)n} y^{kn}}
= \frac{T_1^nx^n}{(xy)^{n-1} y^{kn}}
= \frac{T_1^n}{y^{n(k+1)-1}}  x \, \, \, . $$
Hence it is enough to show that
$\, T_1^n/y^{n(k+1)-1}\,$ is defined on $D(x,y)$, for then we know that
$(x,y)$ is the unit ideal in the ring of global sections $\Gamma(D(x,y),\O_A)$
(and $D(x,y)$ is affine and $A$ is not parasolid).
From $T_1x^{k+1} =z^{n-1}(xy)^k -T_2y^{k+1}$ we get
\begin{eqnarray*}
(T_1x^{k+1})^n &=& (z^{n-1}(xy)^k -T_2y^{k+1})^n
\, = \, y^{kn} (z^{n-1}x^k-T_2y)^n \cr
& =& y^{kn} \sum_{i=0}^n 
\left(
\begin{array}{c} 
n \cr 
i
\end{array}
\right) x^{ki}z^{(n-1)i} (-T_2y)^{n-i}                \cr
&=& y^{kn} \sum_{i=0}^n 
\left(
\begin{array}{c} 
n \cr 
i
\end{array}
\right) x^{ki} z^{n(i-1)} z^{n-i} (-T_2y)^{n-i}  \cr
&=& y^{kn} \sum_{i=0}^n 
\left(
\begin{array}{c} 
n \cr 
i
\end{array}
\right) x^{ki} (xy)^{i-1} (-T_2zy)^{n-i}   \cr
&=& \frac{y^{kn} y^{n-1}}{x} \sum_{i=0}^n
\left(
\begin{array}{c} 
n \cr 
i
\end{array}
\right)
x^{ki+i} (-T_2z)^{n-i} 
\, = \,  \frac{y^{(k+1)n-1}}{x} (x^{k+1} -T_2z)^n \, 
\end{eqnarray*}
and hence
$ T_1^n/y^{n(k+1)-1} = (x^{k+1}-T_2z)^n/x^{n(k+1)+1}$,
so that this function is defined.
\end{example}

\medskip
In the preceeding example it would have been enough to consider only
$k=0$, as the following lemma shows.

\begin{lemma}
\label{pararegulartwo}
Let $R$ denote a local two-dimensional Gorenstein ring.
Let $c =s/xy \in H^2_\fom(R)$ be an element such that
$\Ann \, c= \fom$.
Then $R$ is pararegular if and only if
$D(x,y) \subseteq \Spec \, R[T_1,T_2]/(xT_1+yT_2+s)$ is affine.
\end{lemma}
\proof
We may assume that $R$ is complete.
Due to \ref{pararegularchar} it is enough to show that
if $c=0$ in an $R$-algebra $B$, then $B$ is not universally
parasolid, and for that it is enough to show that
$D(x,y) \subseteq \Spec\, B$ is affine.
If $c=0$ holds in $B$ then this algebra factors through
a forcing algebra for elements $x^{i+1},y^{j+1}; sx^{i}y^{j}$ for some
$i,j \geq 0$.
We have a ring homomorphism
$$R[S_1,S_2]/(sx^{i}y^{j}+x^{i+1}S_1+ y^{j+1}S_2) \lra
R[T_1,T_2]/(s+xT_1+yT_2)$$
given by
$S_1 \mapsto y^{j}T_1$, $S_2 \mapsto x^{i}T_2$.
This mapping is an isomorphism on $D(x,y)$ (check locally on $D(x)$ and $D(y)$),
hence if $D(x,y)$ is affine on the right, it is affine on the left
(In fact $D(x,y) \subseteq  \Spec\, R[T_1,T_2]/(xT_1+yT_2+s)$ is an
affine-linear bundle of rank one over $D(x,y) \subseteq \Spec \, R$
and it is a geometric representation (as a torseur) of the class
$c=s/xy \in H^1(D(x,y), \O_R) \cong H^2_\fom (R)$).
\qed

\section{Expressing datas over $\ZZ$}
\label{expresssec}

Suppose that a certain $R$-algebra is not parasolid.
Then we can express this fact by finitely many
equations over $\ZZ$ and from this we get a lot of other
rings $R'$ where the corresponding algebra is also not parasolid.
This gives in particular relations between the situation in
positive characteristic, zero characteristic and
mixed characteristic.
The following lemma is the most general version of this observation.

\begin{lemma}
\label{triviallemma}
Let $D$ denote a commutative ring and let
$B$ be a finitely generated $D$-algebra.
Suppose that there exists a local Noetherian $D$-algebra $R$
of dimension $d$
such that $B \otimes_D R$ is not parasolid.
Then there exists a finitely generated
$D$-algebra $S=D[X_1, \ldots ,X_d, Y_1, \ldots ,Y_m]/\foa$
such that

\renewcommand{\labelenumi}{(\roman{enumi})}
\begin{enumerate}

\item
$S \subseteq R$ and $X_1, \ldots ,X_d$ are parameters in $R$.

\item
The cohomology class
$1/X_1 \cdots X_d \in H^d_{(X_1, \ldots ,X_d)}(S)$
is not zero.

\item
The cohomology class
$1/X_1 \cdots X_d \in H^d_{(X_1, \ldots ,X_d)}(B \otimes_D S)$
is zero.
\end{enumerate}

Moreover,
for every local Noetherian $S$-algebra $R'$ such that $X_1, \ldots, X_d$
become parameters with non-vanishing paraclass,
the algebra $B \otimes_DR'$ is not parasolid.
\end{lemma}

\proof
Let $B$ be given by $B=D[T_1, \ldots , T_n]/(P_1, \ldots, P_m)$.
Since $B \otimes_D R$ is not parasolid,
there exists parameters
$x_1, \ldots ,x_d$ in $R$ (with non vanishing paraclass) such that
their paraclass does vanish in
$B \otimes _D R=R[T_1, \ldots ,T_n]/(P_1, \ldots , P_m)$.
This means that
$(x_1 \cdots x_d)^k \in (x_1^{k+1}, \ldots ,x_d^{k+1})$
or that we have an equation
$$ (x_1 \cdots x_d)^k
=G_1 x_1^{k+1} + \ldots +G_d x_d^{k+1} +H_1P_1 + \ldots +H_mP_m $$
where $G_i, H_j \in R[T_1, \ldots ,T_n]$.
Let $y_j \in R$ be all the coefficients of
$G_i$ and $H_j$ together.
Now take the finitely generated
$D$-subalgebra $S:=D[x_1, \ldots ,x_d, y_j] \subseteq R$.

(i) is clear by definition and (ii) is clear,
since the paraclass maps to the paraclass in $R$.
(iii). We have the inclusion
$S[T_1, \ldots , T_n] \subseteq R[T_1 , \ldots , T_n]$
and so the last equation holds also in $S[T_1, \ldots , T_n]$,
and hence $(x_1 \cdots x_d)^k \in (x_1^{k+1}, \ldots ,x_d^{k+1})$
holds in $B \otimes _D S =S[T_1, \ldots ,T_n]/(P_1, \ldots , P_m)$.

If $S \ra R'$ is a mapping to a local Noetherian ring
where the images of the $X_i$ are parameters with non-vanishing
paraclass,
then this class vanishes in $B \otimes_D R'$, showing that
this algebra is not parasolid.
\qed

\begin{corollary}
\label{corzet}
Let $D = \ZZ[F_1, \ldots, F_n,H]$ and let
$B= D[T_1, \ldots ,T_n]/(F_1T_1+ \ldots +F_nT_n+H)$ be the universal
forcing algebra.
Let $R$ be a local Noetherian ring and let
$f_1, \ldots, f_n,h \in R$ be elements
so that $R$ is a $D$-algebra via $F_i \mapsto f_i$,
$H \mapsto h$.
Suppose that the forcing algebra $B \otimes_D R$ over $R$
is not parasolid and that the paraclass to
the parameters $x_1, \ldots ,x_d \in R$
is $\neq 0$, but vanishes in the forcing algebra.
Then there exists a finitely generated $\ZZ$-subalgebra $S \subseteq R$,
$$S=\ZZ[ F_1, \ldots, F_n,H,X_1, \ldots ,X_d, Y_j]/\foa$$
such that the cohomology class
$1/X_1 \cdots X_d \in H^d_{(X_1, \ldots ,X_d)}(S)$
is not zero, but it is
zero in
$H^d_{(X_1, \ldots ,X_d)}(B \otimes_D S)$.
\end{corollary}

\proof
This follows from \ref{triviallemma}.
\qed

\begin{remark}
\label{finitelytype}
The previous lemma leads to the following finitely generated
$\ZZ$-algebras (for $d,n,m,k \in \NN$)
$$ S= \ZZ[F_1, \ldots ,F_n, H,X_1, \ldots ,X_d,
C_\nu, C_{1 \nu}, \ldots, C_{d \nu },\, | \nu | \leq m]/ \foa \, , $$
where the relations (the generators of $\foa =\foa_k$) have to be choosen
in such a way that they yield the equation
$$ (X_1 \cdots X_d)^k + X_1^{k+1}\sum_{|\nu| \leq m} C_{1 \nu } T^\nu 
+ \ldots +X_d^{k+1}\sum_{|\nu| \leq m} C_{d \nu } T^\nu
+\sum_{|\nu| \leq m} C_{\nu}T^\nu (\sum_{i=1}^n T_iF_i+H) $$
in the forcing algebra for $F_i, H$ over $S$.
This gives the relations (looking at the coefficients of each monomial
$T^\nu$)
$$ (X_1 \cdots X_d)^k +C_{1 0} X_1^{k+1} + \ldots + C_{d 0}X_d^{k+1}+C_0H  $$
and
$$ C_{1 \nu} X_1^{k+1} + \ldots + C_{d \nu}X_d^{k+1}+C_\nu H 
+F_1 C_{\nu -e_1}+ \ldots +F_n C_{\nu- e_n} \, ,$$
where $|\nu| \leq m$, $e_i$ is the $i$-th unit vector and
$C_{\nu -e_j}=0$ for $\nu_j=0$.
\end{remark}

\medskip
Corollary \ref{corzet} may be applied to show that if some
rings containing a field are pararegular, then also
the corresponding rings in mixed characteristic are pararegular.
This technique may be thought of as a kind of reduction
to the field case.

\begin{proposition}
Let $G \in \ZZ[X,Y,Z]$ be a polynomial
of type $G=XG_1+YG_2+Z^n,\, n \geq 1$
and suppose that
$R=\QQ[X,Y,Z]/(G)$ is pararegular.
Let $C$ denote a three-dimensional regular ring
with maximal ideal $\fom =(x,y,z)$.
Then $R'=C/(G)$ is pararegular except for finitely many residue
characteristics.
\end{proposition}
\proof
$R'$ is a local Gorenstein ring and
the elements $x$ and $y$ are parameters.
The annihilator of the class
$z^{n-1}/xy$ is the maximal ideal $(x,y,z)$.
Due to \ref{pararegulartwo}
we have to show that a geometric realization of this class
over $D(x,y)$ is affine.
Let $B=\ZZ[X,Y,Z][T_1,T_2]/(G,XT_1+YT_2+Z^{n-1})$
be the forcing algebra over $\ZZ$.

Due to our assumption the open subset
$D(X,Y) \subseteq \Spec \, B_\QQ$ is affine.
This means
that there exists an equation
$$ (XY)^k+ C_1X^{k+1} +C_2Y^{k+1} + C_3( XT_1+YT_2+Z^{n-1})=0$$
in $B_\QQ$.
Taking into account the coefficients in $C_i$, this equation holds also
in $B_g$, $0 \neq g \in \ZZ$, hence
$$D(X,Y) \subseteq \Spec \ZZ_g[X,Y,Z][T_1,T_2]/(G, XT_1+YT_2+Z^{n-1})$$
is affine.
If the residue characteristic of $R'$ does not divide $g$,
then we have a factorization $\ZZ_g[X,Y,Z]/(G) \ra R'$
and $D(x,y)$ is affine in the forcing algebra over $R'$.
\qed

\begin{remark}
The previous proposition may be applied
to the equations
$xy +z^{n}=0$, $x^2y+ y^k+ z^2=0$,
$x^2+y^3+z^4=0$, $ x^3+xy^3+y^2+z^2=0$, $x^2+y^3+z^5=0$.
These equations define over $\QQ$ a quotient singularity,
hence they are pararegular. Hence we know by the proposition
that these equations
give also in mixed characteristic pararegular rings except for finitely
many characteristics.
If we want to have more information about the number
$g \in \ZZ$ we need to compute an explicit representation of the unit,
as in example \ref{toricexample}.

This can be obtained in the following way.
We have an explicit quotient map
$$\varphi: \AA^2_\QQ \ra \Spec\, R,\, \, \,
\varphi(u,v)=
(\varphi_1(u,v), \varphi_2(u,v),\varphi_3(u,v)) =(x,y,z)\, .$$
When we pull back the forcing algebra $A$ for $x,y,z^{n-1}$ to the regular ring
$\QQ [u,v]$
we get a forcing algebra $B$ over a regular ring
which has an easy structure,
in particular we can find
functions $p', q'$ defined on
$D(u,v) \subseteq \Spec\, B$
such that $p'u+q'v=1$.
The group acts also on $\Spec \, B$ and on the ring of global sections
$\Gamma(D(u,v),\O_B)$
and so we find also
invariant functions $p,q$ such that $px+qy=1$.
This gives the equation in $A$ showing that $D(x,y)$ is affine.

This method works in principle also in higher dimensions, but then it is not
clear that we only have to look at one forcing algebra.
\end{remark}

\medskip
Another feature of corollary \ref{corzet} is that it allows us
to apply Hochster's finiteness theorem
(see \cite[Theorem 8.4.1]{brunsherzogrev})
to parasolid closure and so to reduce statements in characteristic
zero to positive characteristic.
With this method we may prove
the Theorem of Brian\c{c}on-Skoda in characteristic zero
for parasolid closure.

\begin{proposition}
Let $R$ denote a Noetherian ring containing a field $K$ and let
$I=(f_1, \ldots ,f_n)$ be an ideal.
Then $\overline{I^{n+ w}} \subseteq (I^{w +1})^\pasoclo$.
\end{proposition}
\proof
We may assume that $R$ is local.
In positive characteristic this is a standard result for tight closure
(see \cite[Theorem 5.7]{hunekeapplication}),
so by \ref{tightpasoclop} and \ref{tightpasoclop2}
it is also true for parasolid closure.
Lets assume
that the characteristic of $K$ is zero and assume that
we have a counter example in the local ring $R \supseteq \QQ$, i.e.
$h \in \overline{I^{n+ w}}$ but $h \not\in (I^{w +1})^\pasoclo$.
This means in particular
that the forcing algebra for $f_1, \ldots, f_n,h$ is not
parasolid over $R \ra R'$, but since the condition on the integral closure is
persistent we may assume that $R=R'$. Let $x_1, \ldots ,x_d \in R$
be parameters in $R$ such that their paraclass vanishes
in the forcing algebra.

We have to express the relevant datas in a finitely generated
$\ZZ$-algebra.
The condition $h \in \overline{I^{n+ w}}$ means that
we have an equation
$$h^s + a_1h^{s-1} + \ldots +a_s =0, \, \mbox{ where }
a_i \in (I^{n+w})^i =I^{i(n + w)} \, .$$
This means that we may write
$a_i = \sum_{\nu, \, |\nu|=i(n+w)}
r_\nu f_1^{\nu_1} \cdots f_n^{\nu_n}$.
Furthermore the vanishing of the paraclass has to be expressed as in
\ref{finitelytype}.
Hence everything may be expressed by the $\ZZ$-algebra
$$\ZZ[F_1, \ldots ,F_n,H, R_\nu \,( |\nu|=i(n+w),\,i=1, \ldots ,s),
X_1, \ldots ,X_d, C_{1\nu }, \ldots ,C_{d \nu}, C_\nu]/\foa \, ,$$
and since this algebra has a solution of height $d$ (with respect to
$X_1, \ldots ,X_d$) in a local Noetherian ring over
a field of characteristic zero, it has also a solution of height $d$
over a field of positive characteristic,
see \cite[Theorem 8.4.1]{brunsherzogrev}.
This yields a contradiction.
\qed

\begin{remark}
\label{coloncaptremark}
We take a look at the situation of colon-capturing, i.e. suppose that
$hy_t \in (y_1, \ldots ,y_{t-1})$, where $y_1, \ldots ,y_t$
are parameters. We would like to prove that
$h \in (y_1, \ldots ,y_{t-1})^\pasoclo$.
This is true again via tight closure in positive characteristic, but applying
the method of the last proposition shows only
that the corresponding forcing algebra is parasolid, whereas
universally parasolid is not clear due to the fact that
the condition in colon-capturing is not persistent.
Nethertheless we will show that colon-capturing holds for
parasolid closure in equal characteristic zero in the next section.
\end{remark}

\medskip
Our definition of parasolid closure refers to all local Noetherian rings
in order to make it persistent.
The following proposition describes an important situation where it is enough
to consider only algebras essentially of finite type.

\begin{proposition}
\label{essential}
Let $K$ be a field and let $R$ denote a $K$-algebra of finite type.
Let $A$ be an $R$-algebra of finite type and suppose that
$A'=A \otimes _R R'$ is parasolid for every $R \ra R'$,
where $R'$ is local and essentially of finite type over $K$.
Then $A$ is universally parasolid, i.e. $A'$ is parasolid
for every Noetherian local ring $R'$.
\end{proposition}
\proof
Suppose that
$R \ra R'$
is a Noetherian local ring of dimension $d$ such that $A \otimes_R R'$
is not parasolid. Due to \ref{triviallemma}
there is a finitely generated $R$-subalgebra $S \subseteq R'$
containing elements $x_1, \ldots , x_d \in S$ which are parameters in $R'$
and which have the properties of the lemma.

$S$ is also finitely generated over the field $K$
and the superheight of $(x_1, \ldots ,x_d)$
is exactly $d$. We know due to a
theorem of Koh \cite{kohsuper} that the superheight
equals the finite superheight, hence
there exists a morphism essentially of finite type
$S \ra R''$
($R''$ may be choosen to be the localization at a maximal ideal
of $S'$, $S \ra S'$ of finite type)
such that $(x_1, \ldots ,x_d)$ are parameters in the
local ring $R''$.
But then their paraclass does vanish in $H^d_{\fom''} (A'')$ and we have
a contradiction.
\qed

\begin{corollary}
\label{esspara}
Let $K$ be a field and let $R$ denote a $K$-algebra of finite type.
Let $I \subset R$ be an ideal and let $h \in R$ be an element, and let
$A$ denote the forcing algebra for $h,I$.
If $A'$ is parasolid for every $R \ra R'$ essentially of finite type over $K$,
$R'$ local, then $h \in I^\pasoclo$.
\end{corollary}

\proof
This follows directly from \ref{essential}.
\qed

\section{Generic properties}
\label{generic}

We want to show that tight closure in equal characteristic zero
is contained inside the parasolid closure.
For this we have to consider a relative situation
over a base scheme $\Spec\, D$
and we have to compare the parasolid closure in the generic
fiber with the parasolid closure in the special fibers.
For a $D$-algebra $S$ and a point $\defop \in \Spec\, D$ we denote
by $S_\defop =S \otimes_D \kappa(\defop)$ the ring of the
fiber over $\defop$.
Later on $D$ will be a finitely generated $\ZZ$-domain (often $D=\ZZ$)
so that the closed points have positive characteristic while
the generic point has characteristic zero.

\begin{lemma}
\label{genericlemma}
Let $D$ be a Noetherian domain
with generic point $\eta \in \Spec \, D$.
Let $S$ be an domain of finite type over $D$
such that the fibers are equidimensional of dimension $d$
over a non-empty open subset.
Let $R=S_\eta$ be the ring of the generic fiber,
which is a $d$-dimensional domain of finite type over $Q(D)$.
Let $\fom $ denote a maximal ideal of $R$ and suppose that
$x_1, \ldots ,x_d \in R$ are parameters in $R_\fom$.
Let $A$ be an $R$-algebra of finite type and suppose
that the paraclass $1/x_1 \cdots x_d$ vanishes in $H^d_\fom (A)$.

Then there exists a non-empty
open affine subset $D(g) \subseteq \Spec \, D$,
a prime ideal $\fon \subseteq S$ and an $S$-algebra $B$ of finite type
such that the following hold.

\renewcommand{\labelenumi}{(\roman{enumi})}
\begin{enumerate}

\item
The elements $x_1, \ldots , x_d$ are defined on $S_g$
and they are primary to $\fon \subseteq S_g$.

\item
For every point
$\defop \in D(g)$ the intersection
$\Spec \, S_\defop \cap V(\fon)$ is non-empty and zero-dimensional,
and for every closed point $P \in \Spec \, S_\defop \cap V(\fon)$
the $x_1, \ldots ,x_d$ are parameters in $(S_\defop)_P$.

\item 
$A=B_\eta$.

\item
For every $\defop \in D(g)$
and every closed point $P \in \Spec\, S_\defop \cap V(\fon)$
the paraclass
$1/x_1 \cdots x_d \in H^d_P ((S_\defop)_P)$
vanishes in $H^d_P((B_\defop)_P)$.

\end{enumerate}
\end{lemma}

\proof
During the proof the open subset $D(g)$ gets smaller and smaller.
Of course there exists $0 \neq g \in D$ such that
$x_1, \ldots, x_d \in S_g$.
There exists a unique prime ideal $\fon \subset S$ such that 
$\fon S_\eta = \fom$,
$\fon \cap D = \{ 0 \}$.
We have $(x_1, \ldots, x_d) \subseteq \fon$ in $S_g$ and
$\fon $ is a minimal prime over $(x_i)$.
Let $\fon_i$ denote the other minimal primes over $(x_1, \ldots, x_d)$.
Since these do not survive in $S_\eta$, there exist elements
$ 0 \neq g_i \in D$ and $g_i \in \fon_i$. Making $g_i$ to units,
we may assume that $(x_1, \ldots ,x_d) \subseteq \fon$ is the only
minimial prime.
This gives (i).

The mapping $V(\fon) \ra \Spec \, D$ is finite in the generic point,
hence it is finite over an open non-empty subset,
and $V(\fon) \cap \Spec\, S_\defop$ is zero-dimensional
and non-empty.
If $P$ is a closed point in this intersection, then the $x_i$ generate
up to radical the maximal ideal
$\fom_P$ in $(S_\defop)_P$.
Since the dimension of every component of the fiber
$\Spec\, S_\defop$ is $d$, the height of $P$ in the fiber is $d$,
hence the $x_i$ are parameters.

The construction of $B$ is clear.

Generically we have $(x_1 \cdots x_d)^k \in (x_1^{k+1} , \ldots , x_d^{k+1})$
in $A=B_\eta$. The coefficients of this equation belong to $B_g$
and hence this equation holds over an open neighbourhood.
Then it holds for every $B_\defop$ over $D(g)$ and this shows that
the paraclass $1/x_1 \cdots x_d \in H^d_P ((B_\defop)_P)$
vanishes.
\qed

\begin{remark}
If $D$ is a Hilbert ring so that the closed points
in $\Spec\, D$ are dense, then the condition that the fibers
are equidimensional over an open non-empty subset
is superfluous, since this follows from
\cite[Theorem 14.8]{eisenbud} and the condition that $S$ is a domain.
\end{remark}

\begin{corollary}
\label{parasolidgeneric}
Let $D$ be a Noetherian domain
with generic point $\eta \in \Spec \, D$.
Let $S$ be a domain of finite type over $D$
such that the fibers are equidimensional.
Let $B$ be an $S$-algebra. Suppose that
the set
$\{ \defop \in \Spec \, D: \, S_\defop \ra B_\defop \,
\mbox{ is parasolid }\, \}$
is dense in $\Spec \, D$.
Then $S_\eta \ra B_\eta $ is parasolid in the generic point.
\end{corollary}

\proof
Suppose that $B_\eta$ is not parasolid in the generic point.
This means that there exists a maximal ideal $\fom \subseteq S_\eta=R$
and parameters $x_1, \ldots ,x_d \in R_\fom$ such that their
paraclass vanish in $(B_\eta )_{\fom}$.
We find $0 \neq s \in S$, $s \not\in \fom$
such that $x_1, \ldots ,x_d \in S_s=:S'$.
We may then work with $S'$ instead of $S$.
Due to \ref{genericlemma} the $x_i$ are parameters
in the closed points of the fibers over
an open non-empty subset and their paraclass vanishes in
$B_\defop$, hence $B_\defop $ is not parasolid over $S_\defop$
for an open non-empty subset, yielding a contradiction. 
\qed

\medskip
We obtain an even stronger result if the special points have
positive characteristic.

\begin{corollary}
Let $D$ be a $\ZZ$-domain of finite type
and let $S$ be a domain of finite type over $D$.
Suppose that the fibers $\Spec \, S_\defop$ are normal
{\rm (}or analytically irreducible{\rm )}.
Let $B$ be an $S$-algebra. Suppose that
the set
$\{ \defop \in \Spec \, D: \, S_\defop \ra B_\defop \,
\mbox{ is parasolid}\, \}$
is dense in $\Spec \, D$.
Then $S_\eta \ra B_\eta $ is universally parasolid in the generic point
$\eta \in \Spec\, D$.
\end{corollary}
\proof
First note that the residue class fields of the closed points
in $\Spec\, D$ have positive characteristic.
If for a closed point $\defop$ the $S_\defop$-algebra $B_\defop$
is parasolid, then it is also universally parasolid due to the
assumptions on the fibers and \ref{parasolidpos2} and \ref{parasolidcompl}.
We have to show that
$A=B_\eta$ is universally parasolid over
$R=S_\eta$ and due to
\ref{essential}
it is enough to show that $A'$ is parasolid over
$R'$, where $R'$ is of finite type over $R$.
There exists also an $S$-algebra $S'$ of finite type such that
$R'= S'_\eta$.
Then the $S'$-algebra $B'=B \otimes_S S'$ is again parasolid
over a dense subset of $\Spec\, D$
and then also the generic algebra is parasolid.
\qed

\medskip
We apply these results to forcing algebras and to
the parasolid closure of an ideal.

\begin{proposition}
\label{genericpara}
Let $D$ denote a Noetherian domain
with generic point $\eta$.
Let $S$ be a $D$-domain of finite type
such that the fibers are equidimensional,
and let
$R = S_\eta$ be the ring of the generic fiber.
Let $I \subseteq S$ be an ideal and $h \in S$ and let $B$
be the forcing algebra for $h$ and $I$.
Then the following hold.

\renewcommand{\labelenumi}{(\roman{enumi})}
\begin{enumerate}

\item
If the forcing algebra $A=B_\eta$ is not parasolid in the generic point,
then there
exist a non-empty open subset $D(g) \subseteq \Spec \, D$
such that
the forcing algebra $B_\defop$ is not parasolid
over $S_\defop$ for $\defop \in D(g)$.

\item
Suppose that $h \not\in I^\pasoclo$.
Then there exists a non-empty open subset $D(g) \subseteq \Spec \, D$
such that
the forcing algebra $B_\defop $ is not universally parasolid
over  $S_\defop$, $\defop \in D(g)$
{\rm (}If the characteristic of $\kappa(\defop)$ is positive,
this means that $h \not\in I^\pasoclo$ in $S_\defop$.{\rm)}
\end{enumerate}
\end{proposition}

\proof
The first statement follows from \ref{parasolidgeneric}, so consider (ii).
Since $h \not\in I^\pasoclo$,
the forcing algebra $B_\eta$ is not universally parasolid, hence
we know due to \ref{esspara}
that there exists a ring homomorphism
$R \ra R'$ of finite type to a domain $R'$ and a maximal ideal
$\fom' \subseteq R'$ such that $A'$ is not parasolid over $R'_{\fom'}$.
There exists an $S$-domain $S'$ of finite type such that
$R' =S'_\eta$.
Part (i) shows that there exists a non-empty open subset
$D(g) \subseteq \Spec \, D$
such that $B_\defop'$ is not parasolid over $S_\defop'$ for $\defop \in D(g)$,
hence $B_\defop$ is not universally parasolid over $S_\defop$
for $\defop \in D(g)$.
\qed

\bigskip
There are different ways to define tight closure
in equal characteristic zero all refering to positive characteristic.
We cite the definition from \cite[Appendix, definition 3.1]{hunekeapplication}.

\begin{definition}
\label{tightdefinition}
Let $C$ be a locally excellent Noetherian
ring containing $\QQ$, let $h \in C$ and let
$I=(f_1, \ldots ,f_n) \subseteq C$ be an ideal.
Then $h$ belongs to the tight closure of $I$, $h \in I^*$,
if there exists a finitely generated $\ZZ$-subalgebra $S \subseteq C$
with $h \in S$ and $f_1, \ldots ,f_n \in S$
such that $h_\kappa \in I_\kappa^*$ in $S_\defop$
for almost all $\kappa = \kappa(\defop)$, $\defop \in \Spec\, \ZZ$.
\end{definition}

\begin{proposition}
\label{tightgeneric}
Let $\ZZ \subseteq S$ be a $\ZZ$-domain of finite type
and let $I \subseteq S$ be an ideal, $h \in S$.
Set
$R=S_\QQ$ and suppose that $h \not\in I^\pasoclo$ in $R$.
Then there exists an open non-empty subset $U \subseteq \Spec \, \ZZ$ such that
$h \not\in I^*$ for every closed point $\defop \in U$.
\end{proposition}

\proof
This follows from \ref{genericpara} and \ref{tightpasoclop}.
\qed

\begin{remark}
This result means that if $h \in I_{\kappa (\defop)}^*$ holds for a dense
(=infinite) subset of points $\defop $ with positive characteristic,
then $h \in I_{\QQ}^\pasoclo$ in the generic point (of characteristic zero).
This is in particular true if $h \in I_{\kappa (\defop)}^*$
holds for an open non-empty subset \--- and this is the definition of tight closure
in characteristic zero. It is not clear whether
$h \in I_{\kappa (\defop)}^*$ can be true for a dense set without being true
for almost all $\kappa(\defop)$ (see \cite[Appendix 1.8]{hunekeapplication}).
\end{remark}

\begin{corollary}
\label{tightparazero}
Let $K$ denote a field of characteristic zero and let
$C$ be a Noetherian locally excellent $K$-algebra
and let $I \subseteq C$ be an ideal.
Then $I^* \subseteq I^\pasoclo $.
\end{corollary}

\proof
Suppose that $h \in I^*$.
We may express the relevant datas in a finitely
generated $\ZZ$-algebra $S \subset C$ such that
$S_\QQ \subseteq C$ and where $h_\kappa \in I_\kappa^*$
holds for almost all $\kappa$.
We may assume that $C$ and $S$ are domains.
Then due to \ref{tightgeneric} we have $h \in I^\pasoclo$ in
$S_\QQ$ and then due to the persistence of parasolid closure
this holds also in $C$.
\qed

\begin{remark}
It follows from \ref{tightparazero} that
also colon capturing (and again the Theorem of Brian\c{c}on-Skoda)
holds in equal characteristic zero for parasolid closure.
Of course one would like to prove these facts also without
referring to positive characteristic, and also for mixed characteristic.
\end{remark}

%===========================================================

\end{document}